\documentclass[reqno,11pt]{amsart}

\usepackage{amsfonts,amsmath,amssymb,fullpage,graphicx,color}
\usepackage{verbatim}

\newcommand{\de}{\mathrm{d}}

\newcommand{\nn}{\nonumber}
\newcommand{\ug}{\!\!\!\!&=&\!\!\!\!}

\begin{document}

\makeatletter
\title{On the logarithm of the derivative operator}
\author{D. Babusci}
\address{INFN - Laboratori Nazionali di Frascati, v. le E. Fermi, 40, IT 00044 Frascati (Roma), Italy}
\email{danilo.babusci@lnf.infn.it}
\author{G. Dattoli}
\address{ENEA - Centro Ricerche Frascati, v. le E. Fermi, 45, IT 00044 Frascati (Roma), Italy;  
Universit\'e Paris XIII, LIPN, Institut Galil\'ee, CNRS UMR 7030, 99 Av. J.-B. Clement, F 93430 
Villetaneuse, France}
\email{giuseppe.dattoli@enea.it}

\subjclass[2000]{Primary 33, Secondary 26A33}
\keywords{Operational calculus, Integrals, Bessel functions}

\begin{abstract}
We study the properties of the logarithm of the derivative operator and show that its action on a constant is not zero, 
but yields the sum of the logarithmic function and the Euler-Mascheroni constant. We discuss more general aspects 
concerning  the logarithm of an operator for the study of the properties of the Bessel functions.
\end{abstract}

\maketitle

In this note we deal with the operator 
\begin{equation}
\label{eq:Loper}
\hat{L}_\partial = \ln \partial_x
\end{equation}
and establish how it acts on unity, on monomials and, more in general, on a function.

We define the logarithm according to the identity
\begin{equation}
\label{eq:Logx}
\ln x = \lim_{\nu \to 0} \frac{x^\nu - 1}{\nu} \qquad\qquad (x > 0,\;\nu > 0)\,, \nn
\end{equation}
and, thus
\begin{equation}
\label{eq:Log1}
\hat{L}_\partial\,1 = \lim_{\nu \to 0}\,\left(\frac{\partial_x^\nu - 1}{\nu}\right)\,1\,.
\end{equation}
Before taking the limit, we treat $\nu$ as a generic real number, and write \cite{OldSpa}
\begin{equation}
\label{eq:dfrac}
\partial_x^\nu\,x^\mu = \frac{\Gamma (\mu + 1)}{\Gamma (\mu - \nu + 1)}\,x^{\mu - \nu}
\end{equation}
that substituted in eq. \eqref{eq:Log1} ($\mu = 0$) gives
\begin{equation}
\hat{L}_\partial\,1 = \lim_{\nu \to 0}\,\frac1{\nu}\left[\frac{x^{- \nu}}{\Gamma (1 - \nu)} - 1\right]\,.
\end{equation}
By applying the l'H\^opital rule, we get
\begin{equation}
\label{eq:LimLog}
\hat{L}_\partial\,1 = \lim_{\nu \to 0}\,x^{- \nu}\left[\frac{- \Gamma (1 - \nu)\,\ln x + \Gamma^\prime (1 - \nu)}{\Gamma^2 (1- \nu)}\right] 
= - \ln x - \gamma
\end{equation}
where $\gamma$ is the Euler-Mascheroni constant.

The same procedure can be applied to the study of the action of the operator \eqref{eq:Loper} on monomials. For $n$ integer greater 
than 1, we get 
\begin{eqnarray}
\label{eq:Lopxn}
\hat{L}_\partial\,x^n \ug \lim_{\nu \to 0}\,\frac1{\nu}\left[\frac{n!\,x^{n - \nu}}{\Gamma (n - \nu + 1)} - 1\right] \nn \\
\ug \lim_{\nu \to 0}\,x^{n - \nu}\left[\frac{- n!\,\Gamma (n - \nu + 1)\,\ln x + \Gamma^\prime (n - \nu + 1)}{\Gamma^2 (n - \nu + 1)}\right] \nn \\
\ug x^n\,\left[\psi (n + 1) - \ln x\right] 
\end{eqnarray}
where, according to the definition of the digamma function $\psi (z)$ \cite{AbrSte}
\begin{equation}
\psi (n + 1) = \frac1{n!}\,\Gamma^\prime (n + 1) = - \gamma + h_n \qquad \qquad h_n  = \sum_{k = 1}^n \frac1k
\end{equation}

More in general, for a function $f (x)$ that can be expressed as a power series, i.e. 
\begin{equation}
f (x) = \sum_{n = 0}^\infty \frac{a_n}{n!}\,x^n\,,
\end{equation}
by assuming that the operator and the summation commutes, as a consequence of eq. \eqref{eq:Lopxn} we find
\begin{equation}
\hat{L}_\partial\,f(x) = - (\gamma + \ln x)\,f(x) + \sum_{n = 0}^\infty \frac{a_n}{n!}\,h_n\,x^n\,.
\end{equation}
A significant consequence of the previous formula is the following identity
\begin{equation}
[\hat{L}_\partial + \ln (2\,\sqrt{x})]\,I_0 (\sqrt{x}) =  K_0 (\sqrt{x})
\end{equation}
where $I_0$ and $K_0$ are, respectively, the 0-th order modified Bessel function and Macdonald function. 

Just as a final example, let us now calculate $\hat{L}_\partial\,\ln x$. One has
\begin{equation}
\hat{L}_\partial\,\ln x = \lim_{\nu \to 0}\,\left(\frac{\partial_x^\nu\,\ln x - \ln x}{\nu}\right)\,.
\end{equation}
By exploiting eq. \eqref{eq:dfrac}, it is easy to show that (see also Ref. \cite{OldSpa}, p. 103)
\begin{equation}
\partial_x^\mu\,\ln x = - \frac1{x^\mu\,\Gamma (1 - \mu)}\,\left[\gamma - \ln x + \psi (1 - \mu)\right]\,,
\end{equation}
and therefore we obtain
\begin{equation}
\hat{L}_\partial\,\ln x = - \zeta (2) - (\gamma + \ln x)\,\ln x 
\end{equation}
where $\zeta$ is the Riemann zeta function.

It is worth stressing that, according to the previous discussion, the operator $\hat{L}_\partial$ does not commute with the 
derivative itself. For example, as a consequence of eqs.\eqref{eq:LimLog} and \eqref{eq:Lopxn}, one has
\begin{equation}
\left[ \hat{L}_\partial, \partial_x \right]\,x^n = \frac1{x}\,\delta_{n0}\,.
\end{equation}
This is not surprising since the fractional derivative, at least according to the definition given in eq. \eqref{eq:dfrac}, displays 
an analogous behaviour.  It is also interesting to note that the operator $\hat{L}_\partial$ once acting on an entire function produces 
a new function with a logarithmic singularity. The action of its inverse $(\hat{L}_\partial)^{- 1}$, if it exists, on a function with logarithmic 
singularities could, in principle, be a tool to remove these divergences. Further comments on this last point and on the relevant 
applications in the theory of renormalization will be presented elsewhere.
\vspace{0.5cm}

In a recent series of papers \cite{BabDat,Gorska}, the theory of Bessel functions has been reformulated using the following notation
\begin{equation}
\label{eq:bess}
J_\alpha (x) = \left(\hat{c}\,\frac{x}2\right)^\alpha \,\exp\left\{- \hat{c}\,\left(\frac{x}2\right)^2\right\}\,\varphi (0)\,,
\end{equation}
where $\varphi (r) = 1/\Gamma (r + 1)$ and the operator $\hat{c}$ is defined by the following identity
\begin{equation}
\label{eq:umbra}
\hat{c}^{\,\mu}\,\varphi (0) = \varphi (\mu)\,.
\end{equation}
According to this formalism, the derivatives of a Bessel function with respect to its index can be written as
\begin{equation}
\label{eq:derBes}
\partial_\alpha \, J_\alpha (x) =  \left(\hat{c}\,\frac{x}2\right)^\alpha \,\exp\left\{- \hat{c}\,\left(\frac{x}2\right)^2\right\}\,
\left(\ln \frac{x}2 + \ln \hat{c}\right)\,\varphi (0)\,,
\end{equation}
that faces us with the necessity of evaluating the action of the logarithm of the operator $\hat{c}$ on $\varphi (0)$. 
In analogy with eq. \eqref{eq:Log1}, we set
\begin{equation}
(\ln \hat{c})\, \varphi (0) = \lim_{\nu \to 0} \left(\frac{\hat{c}^{\,\nu} - 1}{\nu}\right)\,\varphi (0)
\end{equation}
i.e., taking into account eq. \eqref{eq:umbra}
\begin{equation}
(\ln \hat{c})\, \varphi (0) = \lim_{\nu \to 0} \left[\frac{\varphi (\nu) - \varphi (0)}{\nu}\right]
\end{equation}
that inserted in eq. \eqref{eq:derBes}, yields
\begin{equation}
\partial_\alpha \,J_\alpha (x) = \left(\ln \frac{x}2\right)\,J_\alpha (x) + \sum_{k = 0}^\infty \frac{(- 1)^k}{k!}\,\left(\frac{x}2\right)^{2\,k + \alpha}\,
\lim_{\nu \to 0} \left[\frac{\varphi (k + \nu + \alpha) - \varphi (k + \alpha)}{\nu}\right]
\end{equation}
from which, by keeping the appropriate limit, we obtain
\begin{equation}
\partial_\alpha \,J_\alpha (x) = \left(\ln \frac{x}2\right)\,J_\alpha (x) - \sum_{k = 0}^\infty \frac{(- 1)^k\,\psi (k + \alpha + 1)}{k!\,\Gamma (k + \alpha + 1)}\,
\left(\frac{x}2\right)^{2\,k + \alpha}\,,
\end{equation}
that is a well known result \cite{AbrSte}.

As further example of known identity which can be derived using the same method, we consider the following integral 
\begin{equation}
B_\alpha = \int_0^\infty \de x \,J_\alpha (x)\,\ln x\,.
\end{equation}
Following the definitions \eqref{eq:Logx}, \eqref{eq:bess}, we find
\begin{eqnarray}
B_\alpha \ug \hat{c}^{(\alpha + 1)/2}\,\Gamma (\frac{\alpha + 1}2)\,\left[\ln 2 - \frac12\,\ln \hat{c} + 
\frac12\,\psi (\frac{\alpha + 1}2)\right]\,\varphi (0) \nn \\
\ug \ln 2 +  \psi (\frac{\alpha + 1}2)\,,
\end{eqnarray}
that is the integral \textbf{6.771} of Ref. \cite{GradRyz}. For $\alpha = 0$, one has
\begin{equation}
B_0 = - \gamma - \ln 2\,.
\end{equation}

We believe that the results we have presented in this paper offer a very effective tool for treating a large body of problems 
in pure and applied mathematics. A more systematic analysis will be presented in a forthcoming investigation.

\section*{Acknowledgments}
One of us (G. D.) recognizes the warm hospitality and the financial support of the University Paris XIII, whose stimulating atmosphere provided 
the necessary conditions for the elaboration of the ideas leading to this paper.

\end{document}